\documentclass{amsart}
\usepackage{amsmath}
\usepackage{amscd}
\theoremstyle{definition}
\theoremstyle{remark}
\numberwithin{equation}{section}

\begin{document}
\title[A mixed problem for a hyperbolic equation]{A mixed problem for a
Boussinesq hyperbolic equation with integral condition}
\author{Said Mesloub (1)}
\address{1.king Saoud university,departement of mathematics Riadh,Arabie Saoudite.}
\address{2.Institut Camille Jordan,Université Lyon1,France.}
\email{1. mesloubs@yahoo.com  ,2.  amansour@math.univ-lyon1.fr }
\author{A. Mansour (2)}
\curraddr{}
\date{Octobre 30, 2008}
\maketitle

{\LARGE \begin{abstract}
A hyperbolic problem wich combines a classical(Dirichlet) and a non-local
contraint is considered.The existence and uniqueness of strong solutions are
proved,we use a functionnal analysis method based on a priori estimate and
on the density of the range of the operator generated by the considered
problem.
\end{abstract}
\vspace{1cm}
{\bfseries 1.INTRODUCTION }\\

\vspace{0.3cm}

The first study of evolution problems with a nonlocal condition - the so
called energy specification - goes back to Cannon$\left[ 5\right]$  ,1963
Using an integral condition $,$we proved the existence and uniqueness of the
solution of a mixed problem wich combine a classical ( Dirichlet )and an
integral condition for the equation . Problems involving local and integral
condition for hyperbolic equations are investigated by the energy
inequalities method in $\left[ 1\right] , \left[ 6\right] ,\left[ 7\right] ,%
\left[ 8\right] ,\left[ 9\right] ,\left[ 10\right] ,\left[ 11\right] $,\\
$\left[12\right] .$In this paper ,we prove the existence and uniqueness of the
solution for the mixed problem $\left( 1\right) -\left( 5\right) .$Our proof
is based on a priori estimate and on the fact that the range of the operator
generated by the considered problem is dense.
\newpage

\vspace{1cm}
{\bfseries 2.Formulation of the problem }\\

\vspace{0.3cm}
In the region $Q=\left( 0,l\right) \times \left( 0,T\right) ,$\textbf{\ }%
with $l<\infty $\textbf{\ }and\textbf{\ }$T<\infty ,$\textbf{\ }we shall
consider the problem 

\vspace{0.3cm}
\begin{equation}
Lu=u_{tt}-\left( b(x,t)u_{x}\right) _{x}-\beta \frac{\partial ^{4}u}{%
\partial t^{2}\partial x^{2}}=f(x,t),\forall (x,t)\in Q  \tag{1}
\end{equation}
\begin{equation}
l_{1}u=u\left( x,0\right) =\varphi _{1}\left( x\right) , \hspace*{0.5cm}  x\in \left(
0,l\right)  \tag{2}
\end{equation}
\begin{equation}
l_{2}u=u_{t}\left( x,0\right) =\varphi _{2}\left( x\right) ,  \hspace*{0.5cm} x\in
\left( 0,l\right)  \tag{3}
\end{equation}
\begin{equation}
u\left( 0,t\right) =0 ,\hspace*{1.5cm} t\in \left( 0,T\right)  \tag{4}
\end{equation}
\begin{equation}
\int_{0}^{l}xu\left( x,t\right) dx=0,  \hspace*{0.5cm} t\in \left( 0,T\right)  \tag{5}
\end{equation}

\vspace{0.5cm}
where $\beta $\ $\in IR_{+}^{*}\;$and $b(x,t)$\ and its derivatives satisfy
the conditions:

\vspace{0.1cm}
C$_{1}:b_{0}\leq b(x,t)\leq b_{1}$\ , $b_{t}\left( x,t\right) \leq b_{2}$\ , 
$b_{x}\left( x,t\right) \leq b_{3},$\\ for any $(x,t)\in \overline{Q}$,%
\newline
C$_{2}\ :$\ $b_{tt}\left( x,t\right) \leq b_{4}$\ , $b_{xt}\left( x,t\right)
\leq b_{5}$\ ,for any $(x,t)\in \overline{Q}$.\newline

\vspace{0.1cm}
The functions $f$, $\varphi _{1}$\ and $\varphi _{2}$  are known functions
wich satisfy the compatibility conditions:

\vspace{0.1cm}
$\varphi _{1}\left( 0\right) =\varphi _{2}\left( 0\right)
=\int_{0}^{l}x\varphi _{1}\left( x\right) dx=\int_{0}^{l}x\varphi _{2}\left(
x\right) dx=0$.\newpage

{\bfseries  3. functional Spaces }\newline

\vspace*{0.3cm}
The problem (1)-(5) can be put in the following operator form:\hspace*{0.5cm} $Lu=\mathcal{F}
$, \hspace*{0.1cm} $u\in D(L)$,\hspace*{0.2cm}where:
\begin{center}
$\;Lu=\left( \mathcal{L}u,l_{1}u,l_{2}u\right) $  and  $%
\mathcal{F}=\left( f,\varphi _{1},\varphi _{2}\right) $. 
\end{center}

\vspace*{0.2cm}
The operator$\;L\;$is considered from $B\;$to $H$,where $B\;$ is the Banach
space consisting of functions $u\in L^{2}(Q)$, satisfying conditions
(4)and(5) with the finite norm:
\begin{center}
$$
\left\| u\right\| _{B}^{2}\mathbf{=}%
\sup_{0\leq \tau \leq T} \ \left[ \left\| u(.,\tau )\right\|
_{L^{2}(0,l)}^{2}+\left\| u_{t}(.,\tau )\right\| _{L^{2}(0,l)}^{2}\right]
\nolinebreak 
$$
\end{center}
and $F$\ is the Hilbert space $L^{2}(Q)\times
L^{2}(0,l)\times L^{2}(0,l)$\ equipped with the norm:
\begin{center}
$\left\| \mathcal{F}\right\| _{H}^{2}=\left\| f\right\| _{L^{2}(Q^{\tau
})}^{2}+\left\| \varphi _{1}\right\| _{L^{2}(0,l)}^{2}+\left\| \varphi
_{2}\right\| _{L^{2}(0,l)}^{2}$.
\end{center}
Let $D(L)$\ denote the domain\ of $L$\
,which is the set of all functions $u\in L^{2}(Q)\;$ for which $%
u_{t},u_{x},u_{tx},u_{tt},u_{ttx}\in L^{2}(Q)\;$ and satisfying
conditions(4)and(5).\\

\vspace*{0.5cm}
{\bfseries 4.  A\ priori\ estimate and its consequences }\\

{\bfseries Theorem 1}: For any function $u\in D(L)$ satisfyies conditions C$_{1}$-C$_{2}$%
\ there exists a positive constant $c,$ such that 
\begin{equation}
\left\| u\right\| _{B}\leq c\left\| Lu\right\| _{H},  \tag{7}
\end{equation}
{\bfseries Proof} : We consider the scalair product in $L^{2}(Q^{\tau })$\ of the
operator $\mathcal{L}u\;$and$\;Mu,$  where $Mu=x\Im _{x}^{*}u_{t}-\Im {x}^{*}(\rho u_{t})$,\\
with $Q^{\tau }=(0,l)\times (0,\tau )$ ,\hspace*{0.3cm}$0\leq \tau \leq T$, and $\Im _{x}^{*}v=\int_{x}^{l}v(\xi ,t)d\xi $,
we obtain \newpage

\begin{align}
\left( \mathcal{L}u,Mu\right) _{L^{2}(Q^{\tau })}\tag{8} \newline
& \mathbf{=}\left( u_{tt},x\Im _{x}^{*}u_{t}\right) _{L^{2}(Q^{\tau
})}-\left( \left( b\left( x,t\right) u_{x}\right) _{x},x\Im
_{x}^{*}u_{t}\right) _{L^{2}(Q^{\tau })}  \notag \\
& -\beta \left( u_{ttxx},x\Im _{x}^{*}u_{t}\right) _{L^{2}(Q^{\tau
})}-\left( u_{tt},\Im _{x}^{*}(\rho u_{t})\right) _{L^{2}(Q^{\tau })}\notag\\
& +\left( \left( b\left( x,t\right) u_{x}\right) _{x},\Im _{x}^{*}(\rho
u_{t})\right) _{L^{2}(Q^{\tau })}+\beta \left( u_{ttxx},\Im _{x}^{*}(\rho
u_{t})\right) _{L^{2}(Q^{\tau })}.  \notag
\end{align}

Making use of conditions (2)-(5)and integrating by parts we estabilish the
equalities: 
\begin{gather}
\left( u_{tt},x\Im _{x}^{*}u_{t}\right) _{L^{2}(Q^{\tau })}=\frac{1}{2}%
\left\| \Im _{x}^{*}u_{t}(.,\tau )\right\| _{L^{2}(0,l)}^{2}  \tag{9} \\
-\frac{1}{2}\left\| \Im _{x}^{*}\varphi _{2}\right\|
_{L^{2}(0,l)}^{2}-\left( \Im _{x}^{*}u_{tt},u_{t}\right) _{L_{\rho
}^{2}(Q^{\tau })},  \notag
\end{gather}
\begin{gather}
-\left( \left( b\left( x,t\right) u_{x}\right) _{x},x\Im
_{x}^{*}(u_{t})\right) _{L^{2}(Q^{\tau })}  \tag{10} \\
=\frac{1}{2}\left\| \sqrt{b\left( .,\tau \right) }u(.,\tau )\right\|
_{L^{2}(0,l)}^{2}-\frac{1}{2}\left\| \sqrt{b\left( .,0\right) }\varphi
_{1}\right\| _{L^{2}(0,l)}^{2}  \notag \\
-\frac{1}{2}\left\| \sqrt{b_{t}\left( .,t\right) }u\right\| _{L^{2}(Q^{\tau
})}^{2}-\left( b_{x}\left( x,t\right) u,\Im _{x}^{*}u_{t}\right)
_{L^{2}(Q^{\tau })}  \notag \\
-\left( b\left( x,t\right) u_{x},u_{t}\right) _{L_{\rho }^{2}(Q^{\tau })}, 
\notag
\end{gather}
\begin{gather}
-\beta \left( u_{ttxx},x\Im _{x}^{*}(u_{t})\right) _{L^{2}(Q^{\tau })}=\frac{%
\beta }{2}\left\| u_{t}(.,\tau )\right\| _{L^{2}(0,l)}^{2}  \tag{11} \\
-\frac{\beta }{2}\left\| \varphi _{2}\right\| _{L^{2}(0,l)}^{2}-\beta \left(
u_{ttx},u_{t}\right) _{L_{\rho }^{2}(Q^{\tau })}.  \notag
\end{gather}

\begin{equation}
-\left( u_{tt},\Im _{x}^{*}(\rho u_{t})\right) _{L^{2}(Q^{\tau })}=\left(
\Im _{x}^{*}u_{tt},u_{t}\right) _{L_{\rho }^{2}(Q^{\tau })},  \tag{12}
\end{equation}
\begin{equation}
\left( \left( b\left( x,t\right) u_{x}\right) _{x},\Im _{x}^{*}(\rho
u_{t})\right) _{L^{2}(Q^{\tau })}=\left( b\left( x,t\right)
u_{x},u_{t}\right) _{L_{\rho }^{2}(Q^{\tau })},  \tag{13}
\end{equation}
\begin{equation}
\beta \left( u_{ttxx},\Im _{x}^{*}(\rho u_{t})\right) _{L^{2}(Q^{\tau
})}=\beta \left( u_{ttx},u_{t}\right) _{L_{\rho }^{2}(Q^{\tau })}.  \tag{14}
\end{equation}
Combining equalities (9)-(14) and(8)we obtain : 
\begin{gather}
\frac{1}{2}\left\| \Im _{x}^{*}u_{t}(.,\tau )\right\| _{L^{2}(0,l)}^{2}+%
\frac{1}{2}\left\| \sqrt{b\left( .,t\right) }u(.,\tau )\right\|
_{L^{2}(0,l)}^{2}  \tag{15} \\
\;\;\;\;\;\;\;\;\;\;\;\;\;\;\;+\frac{\beta }{2}\left\| u_{t}(.,\tau
)\right\| _{L^{2}(0,l)}^{2}  \notag \\
=\frac{1}{2}\left\| \Im _{x}^{*}\varphi _{2}\right\| _{L^{2}(0,l)}^{2}+\frac{%
1}{2}\left\| \sqrt{b\left( .,t\right) }\varphi _{1}\right\| _{L^{2}(0,l)}^{2}
\notag \\
+\frac{\beta }{2}\left\| \varphi _{2}\right\| _{L^{2}(0,l)}^{2}+\frac{1}{2}%
\left\| \sqrt{b_{t}}u\right\| _{L^{2}(Q^{\tau })}^{2}+\left( b_{x}\left(
x,t\right) u,\Im _{x}^{*}u_{t}\right) _{L^{2}(Q^{\tau })}  \notag \\
+\left( \mathcal{L}u,x\Im _{x}^{*}u_{t}\right) _{L^{2}(Q^{\tau })}-\left( 
\mathcal{L}u,\Im _{x}^{*}(\rho u_{t})\right) _{L^{2}(Q^{\tau })}\text{.} 
\notag
\end{gather}
By applying the Cauchy inequality to the last three terms on the right-hand
side of the inequality (15) and making use conditions C$_{1},$combining with
(15),we obtain 
\begin{gather}
\left\| u(.,\tau )\right\| _{L^{2}(0,l)}^{2}+\left\| u_{t}(.,\tau )\right\|
_{L^{2}(0,l)}^{2}  \tag{16} \\
\;\;\;\;\;+\left\| \Im _{x}^{*}u_{t}(.,\tau )\right\| _{L^{2}(0,l)}^{2} 
\notag \\
\leq k\left[ \left\| f\right\| _{L^{2}(Q^{\tau })}^{2}+\left\| \varphi
_{1}\right\| _{L^{2}(0,l)}^{2}+\left\| \varphi _{2}\right\|
_{L^{2}(0,l)}^{2}\right.  \notag \\
\left. \left\| u\right\| _{L^{2}(Q^{\tau })}^{2}+\left\| u_{t}\right\|
_{L^{2}(Q^{\tau })}^{2}+\left\| \Im _{x}^{*}u_{t}\right\| _{L^{2}(Q^{\tau
})}^{2}\right] \text{.}  \notag
\end{gather}
where $\;k=\frac{\max \left( 2,b_{1},\beta
+l^{2},b_{3}^{2}+b_{2},l^{4}\right) \;}{\min (1,b_{0},\beta )}\;.$\\
Applying the Gronwall lemma to(16),and elimining the term $\left\| \Im
_{x}^{*}u_{t}(.,\tau )\right\| _{L^{2}(0,l)}^{2}$\ of the left-hand side of
the inequality we obtain 
\begin{gather}
\left\| u(.,\tau )\right\| _{L^{2}(0,l)}^{2}+\left\| u_{t}(.,\tau )\right\|
_{L^{2}(0,l)}^{2}  \tag{17} \\
\leq k\exp (kT)\left( \left\| f\right\| _{L^{2}(Q^{\tau })}^{2}+\left\|
\varphi _{1}\right\| _{L^{2}(0,l)}^{2}+\left\| \varphi _{2}\right\|
_{L^{2}(0,l)}^{2}\right) .  \notag
\end{gather}
\newpage
Since the left-hand side of (17)does not depend on $\tau $,we take the
supremum with $\tau $\ from $0$\ to $T$, then the estimate (7) follows with $%
c=\sqrt{k}\exp (k\frac{T}{2})$.\\

\vspace*{1cm}
{\bfseries 5. Solvability of the problem } \\

{\bfseries Proposition 1. } The operator $L$\ acting from $B$\ to $H$\ have a closure.

{\bfseries Proof. } ( see $\left[ 3\right] $\ )

Let be $\overline{L}$\ the closure of $L$, $D(\overline{L})$\ its domain .%
\newline
{\bfseries Definition }.The solution of $\overline{L}u=F$\ for any $u\in D(\overline{L})$%
\ is strong solution of problem(1)-(5). we take the limit in the inequality
(7) ,we obtain$\left\| u\right\| _{B}\leq c\left\| \overline{L}u\right\|
_{H} $,$\forall u\mathbf{\in }D\mathbf{(}\overline{L}\mathbf{)}$.\newline
From the inequality we have:

{\bfseries Corollary 1}: The strong solution of problem (1)-(5) when it exists, it's
unique, and depends continuly of data $f,\varphi _{1},\varphi _{2}$.

{\bfseries Corollary 2 }:The set of values $R(\overline{L})$\ of the operator $\overline{L%
}$\ is equal to the closure $\overline{R(L)}$\ of $R(L)$.

\vspace*{0.3cm}
{\bfseries Theorem 2}: If the conditions C$_{1}$-C$_{2}$\ are satisfying ,then for any $%
\mathcal{F=}\left( f,\varphi _{1},\varphi _{2}\right) \in H,\;$there exists
a strong unique solution $u=\overline{L}^{-1}\mathcal{F}=\overline{L^{-1}}%
\mathcal{F\;}$of the probleme (1)-(5) where the estimate \ $\left\|
u\right\| _{B}\leq c\left\| \mathcal{F}\right\| _{H}\;$is satisfying ,where $%
c$\ is a positive constant does not depends of $u$.\newline

\vspace*{0.2cm}
{\bfseries Proof}: From(22) we conclude that the operator $\overline{L}$\ acting from $%
D(\overline{L})$\ in $R(\overline{L})$\ have an inverse $\overline{L}^{-1}$,
and from corollary 2, we conclude that the range $R(\overline{L})$\ of the
operator $\overline{L}$\ is closed. Then we will be proove the density of
the set $R(L)$\ in the space $H$\ (i.e) $\overline{R(L)}=H.$

For this we need the following proposition :

{\bfseries Proposition 2}:If,for all functions $u\in D_{0}(L),$ where\\
$$ D_{0}(L)=  \{ u/\text{}u\in D(L): l_{1}u=l_{2}u=0 \}  $$
and for some function $\omega \in L^{2}(Q),$ we have 
\begin{equation}
\left( \mathcal{L}u,\omega \right) _{L^{2}(Q)}=0\text{,}  \tag{18}
\end{equation}
then $\omega \ $  vanishes almost everywhere in $Q$.\newline
{\bfseries Proof of the proposition 2 }: The relation (18) is given for all $u\in
D_{0}(L) $,we can express it in a particular form .Let $u_{tt}$\ be a
solution of : 
\begin{equation}
b\left( \sigma ,t\right) \left[ x\Im _{x}^{*}u_{tt}-\Im _{x}^{*}(\rho u_{tt})%
\right] =h(x,t),  \tag{19}
\end{equation}
where $\sigma $\ is a constant in $\left( 0,l\right) $\ and $%
h(x,t)=\int_{t}^{T}\omega (x,\tau )d\tau $.\newline
And let $u$ be the fonction\ defined by: 
\begin{equation}
u\mathbf{=}\left\{ 
\begin{array}{c}
0\text{,\qquad \qquad \qquad \qquad si }0\leq t\leq s\text{,} \\ 
\int_{s}^{t}(t-\tau )u_{\tau \tau }d\tau \text{,\qquad si }s\leq t\leq T%
\text{.}%
\end{array}
\right.  \tag{20}
\end{equation}
(19) and (20) follows $u$\ is in $D_{0}(L)$\ and: 
\begin{gather}
\omega (x,t)=\Im _{x}^{*^{-1}}h  \tag{21} \\
=-\left[ b\left( \sigma ,t\right) \left( x\Im _{x}^{*}u_{tt}-\Im
_{x}^{*}(\rho u_{tt})\right) \right] _{t}  \notag \\
=\left[ b\left( \sigma ,t\right) \Im _{x}^{*}(\rho -x)u_{tt}\right] _{t}. 
\notag
\end{gather}

To continue the proof we need the following lemma :\newline
{\bfseries Lemma 2}. The function $\omega $\ defined by (21),belongs to the space $%
L^{2}(Q)$.\newline
{\bfseries Proof of lemma 2 }:We start with the proof of this inequality $\left\Vert \Im _{x}^{\ast }(\rho -x)u_{tt}\right\Vert _{L^{2}(0,l)}^{2}\mathbf{\leq }\frac{l^{4}}{12}\left\Vert u_{tt}\right\Vert _{L^{2}(0,1)}^{2}%
\mathbf{.}$\\

\vspace*{0.3cm}
From this inequality and since the conditions C$_{1}$\ are satisfied we
conclude that $\;b_{t}\left( \sigma ,t\right) \Im _{x}^{\ast }(\rho -x)u_{tt}$
belongs to $L^{2}(Q)\mathbf{.}$\nolinebreak \\

\vspace*{0.3cm}
Because  $\omega (x,t)=\left[b\left(\sigma ,t\right)\Im _{x}^{\ast }(\rho -x)u_{tt}\right]_{t}=$%
\newline
$=b_{t}\left( \sigma ,t\right) \Im _{x}^{\ast }(\rho -x)u_{tt}+b\left(
\sigma ,t\right) \Im _{x}^{\ast }(\rho -x)u_{ttt},$ \hspace*{0.2cm} then we will be prooved
that:\hspace*{0.2cm}$b\left( \sigma ,t\right) \Im _{x}^{\ast }(\rho -x)u_{ttt}\in L^{2}(Q).$ \\

\vspace*{0.3cm}
For this we introduce the $t$-averaging op\'{e}rators $\rho _{\varepsilon }$ of the form
\begin{center}
 $(\rho _{\varepsilon }f)(x,t)=\frac{1}{\varepsilon }\int_{0}^{T}\omega (\frac{t-s}{\varepsilon })f(x,s)ds$,
\end{center}
 where$\;\omega \in C_{0}^{\infty }(0,T), \omega \geq 0,$\\
$\int_{-\infty }^{+\infty }\omega(s)ds=1$, \smallskip $\omega \equiv 0$ \ for $t\leq 0$\  and $t\geq T,$ \\
\vspace*{0.3cm}
applying the operators $\rho _{\varepsilon \text{ }}$\ and $\frac{\partial }{\partial
t}$\ to the equation
\begin{center}
$-b\left( \sigma ,t\right) \Im _{x}^{\ast }(\rho -x)u_{tt}=h(x,t)\mathbf{,}$ 
\end{center}
we obtain
\begin{center}
$\frac{\partial }{\partial t}\left( -b\left( \sigma ,t\right) \Im
_{x}^{\ast }(\rho -x)u_{tt}\right) =$\newline
$\frac{\partial }{\partial t}\left[ -b\left( \sigma ,t\right) \Im
_{x}^{\ast }(\rho -x)u_{tt}+\rho _{\varepsilon }\left( b\left( \sigma
,t\right) \Im _{x}^{\ast }(\rho -x)u_{tt}\right) \right] -\frac{\partial }{%
\partial t}\rho _{\varepsilon }h $.
\end{center}
\vspace*{0.3cm}
Then$\;\;\;\;\;\;\;\;\;\;\;\;\;\;\;\;\;\left\Vert b\left( \sigma ,t\right)
\Im _{x}^{\ast }(\rho -x)u_{tt}\right\Vert _{L^{2}(Q)}^{2}$\newline
$\leq 2\left\Vert \frac{\partial }{\partial t}\left[ b\left( \sigma
,t\right) \Im _{x}^{\ast }(\rho -x)u_{tt}-\rho _{\varepsilon }\left( b\left(
\sigma ,t\right) \Im _{x}^{\ast }(\rho -x)u_{tt}\right) \right] \right\Vert
_{L^{2}(Q)}^{2}$

$\;\;\;\;\;\;\;\;\;\;\;\;\;\;\;\;\;\;\;\;\;\;\;\;+2\left\| \frac{\partial }{%
\partial t}\rho _{\varepsilon }h\right\| _{L^{2}(Q)}^{2}$.\\

\vspace*{0.3cm}
Since $\rho _{\varepsilon }f\underset{\varepsilon \rightarrow 0}{%
\longrightarrow }f$, and \ $\frac{\partial }{\partial t}\left( b\left(
\sigma ,t\right) \Im _{x}^{*}(\rho -x)u_{tt}\right) $\ is bounded in $%
L^{2}(Q)$, then $\omega \in L^{2}(Q)$.\\
\vspace*{0.3cm}
Now we return to the 2$^{nd}$ proposition ,we remplace $\omega $\ in (18) by its representation given by
(21) we have: 
\begin{gather}
\left( u_{tt},\left[ b\left( \sigma ,t\right) \Im _{x}^{*}(\rho -x)u_{tt}%
\right] _{t}\right) _{L^{2}(Q)}  \tag{22} \\
=\left( \left( b\left( x,t\right) u_{x}\right) _{x},\left[ b\left( \sigma
,t\right) \Im _{x}^{*}(\rho -x)u_{tt}\right] _{t}\right) _{L^{2}(Q)}+  \notag
\\
+\beta \left( u_{ttxx},\left[ b\left( \sigma ,t\right) \Im _{x}^{*}(\rho
-x)u_{tt}\right] _{t}\right) _{L^{2}(Q)}.  \notag
\end{gather}
\newpage
Making use conditions(3)-(5),and from the particular forme of$\;u$ given by
(19) and (20), the equality (22) can be simplified.For this integrating by
parts each term of the equality on the sup-domain$\;$\ $Q_{s}=\left(
0,l\right) \times \left( s,T\right) $ where $0\leq s\leq T\;$\

\begin{gather}
\left( u_{tt},\left[ b\left( \sigma ,t\right) \Im _{x}^{*}(\rho -x)u_{tt}%
\right] _{t}\right) _{L^{2}(Q)}  \tag{23} \\
=\frac{1}{2}\left\| \sqrt{b\left( \sigma ,s\right) }\Im
_{x}^{*}u_{tt}(.,s)\right\| _{L^{2}(0,l)}^{2}-\frac{1}{2}\left\| \sqrt{%
b_{t}\left( \sigma ,.\right) }\Im _{x}^{*}u_{tt}\right\| _{L^{2}(Q_{s})}^{2},
\notag
\end{gather}
\begin{gather}
\left( \left( b\left( x,t\right) u_{x}\right) _{x},\left[ b\left( \sigma
,t\right) \Im _{x}^{*}(\rho -x)u_{tt}\right] _{t}\right) _{L^{2}(Q)} 
\tag{24} \\
=-\frac{1}{2}\left\| \sqrt{b\left( .,T\right) b\left( \sigma ,T\right) }%
u_{t}(.,T)\right\| _{L^{2}(0,l)}^{2}  \notag \\
+\frac{1}{2}\int_{Q_{s}}\left[ 3b_{t}\left( x,t\right) b\left( \sigma
,t\right) +b\left( x,t\right) b_{t}\left( \sigma ,t\right) \right] \left(
u_{t}\right) ^{2}dxdt  \notag \\
-\int_{0}^{l}b_{t}\left( x,T\right) b\left( \sigma ,T\right) u\left(
x,T\right) u_{t}\left( x,T\right) dx  \notag \\
+\int_{Q_{s}}\left[ b_{tt}\left( x,t\right) b\left( \sigma ,t\right)
+b_{t}\left( x,t\right) b_{t}\left( \sigma ,t\right) \right] uu_{t}dxdt 
\notag \\
+\int_{Q_{s}}\left[ b_{x}\left( x,t\right) u_{t}+b_{xt}\left( x,t\right) u%
\right] b\left( \sigma ,t\right) \Im _{x}^{*}u_{tt}dxdt\text{,}  \notag
\end{gather}
\begin{gather}
\beta \left( u_{ttxx},\left[ b\left( \sigma ,t\right) \Im _{x}^{*}(\rho
-x)u_{tt}\right] _{t}\right) _{L^{2}(Q)}  \tag{25} \\
=\frac{\beta }{2}\left\| \sqrt{b_{t}\left( \sigma ,.\right) }u_{tt}\right\|
_{L^{2}(Q_{s})}^{2}-\frac{\beta }{2}\left\| \sqrt{b\left( \sigma ,s\right) }%
u_{tt}\left( .,s\right) \right\| _{L^{2}(0,l)}^{2}.  \notag
\end{gather}
Substitution of(23)-(25)into (22) gives 
\begin{gather}
\frac{1}{2}\left\| \sqrt{b\left( \sigma ,s\right) }\Im
_{x}^{*}u_{tt}(.,s)\right\| _{L^{2}(0,l)}^{2}  \tag{26} \\
+\frac{1}{2}\left\| \sqrt{b\left( .,T\right) b\left( \sigma ,T\right) }%
u_{t}(.,T)\right\| _{L^{2}(0,l)}^{2}  \notag \\
+\frac{\beta }{2}\left\| \sqrt{b\left( \sigma ,s\right) }u_{tt}\left(
.,s\right) \right\| _{L^{2}(0,l)}^{2}  \notag \\
=\frac{1}{2}\left\| \sqrt{b_{t}\left( \sigma ,s\right) }\Im
_{x}^{*}u_{tt}\right\| _{L^{2}(Q_{s})}^{2}+\frac{\beta }{2}\left\| \sqrt{%
b_{t}\left( \sigma ,.\right) }u_{tt}\right\| _{L^{2}(Q_{s})}^{2}  \notag \\
+\frac{1}{2}\int_{Q_{s}}\left[ 3b_{t}\left( x,t\right) b\left( \sigma
,t\right) +b\left( x,t\right) b_{t}\left( \sigma ,t\right) \right] \left(
u_{t}\right) ^{2}dxdt  \notag \\
-\int_{0}^{l}b_{t}\left( x,T\right) b\left( \sigma ,T\right) u\left(
x,T\right) u_{t}\left( x,T\right) dx  \notag \\
+\int_{Q_{s}}\left[ b_{tt}\left( x,t\right) b\left( \sigma ,t\right)
+b_{t}\left( x,t\right) b_{t}\left( \sigma ,t\right) \right] uu_{t}dxdt 
\notag \\
+\int_{Q_{s}}\left[ b_{x}\left( x,t\right) u_{t}+b_{xt}\left( x,t\right) u%
\right] b\left( \sigma ,t\right) \Im _{x}^{*}u_{tt}dxdt\text{ .}  \notag
\end{gather}
By applying the Cauchy inequality and Cauchy inequality with $\varepsilon $\
to estimate the last three terms on the right-hand side of the inequality
(26) and making use conditions C$_{1}-C_{2}$, combining the estimates
and(26)taking into account that $\varepsilon =\frac{b_{0}^{2}}{2b_{1}^{2}}$
we obtain 
\begin{gather}
\frac{b_{0}}{2}\left[ \left\| \Im _{x}^{*}u_{tt}(.,s)\right\|
_{L^{2}(0,l)}^{2}+\frac{b_{0}}{2}\left\| u_{t}(.,T)\right\|
_{L^{2}(0,l)}^{2}+\right.  \tag{27} \\
\left. +\beta \left\| u_{tt}\left( .,s\right) \right\| _{L^{2}(0,l)}^{2}
\right]  \notag \\
\leq \left( b_{1}^{2}+\frac{b_{2}}{2}\right) \left\| \Im
_{x}^{*}u_{tt}\right\| _{L^{2}(Q_{s})}^{2}+\frac{\beta b_{2}}{2}\left\|
u_{tt}\right\| _{L^{2}(Q_{s})}^{2}+  \notag \\
+\frac{b_{2}^{2}+b_{1}^{2}+4b_{1}b_{2}+b_{3}^{2}}{2}\left\| u_{t}\right\|
_{L^{2}(Q_{s})}^{2}+  \notag \\
+\frac{b_{2}^{2}+b_{4}^{2}+b_{5}^{2}}{2}\left\| u\right\|
_{L^{2}(Q_{s})}^{2}+\frac{b_{1}^{2}b_{2}^{2}}{b_{0}^{2}}\left\| u\left(
.,T\right) \right\| _{L^{2}(0,l)}^{2}  \notag
\end{gather}

By virtue of the elementary inequality 
\begin{equation}
\frac{b_{1}^{2}b_{2}^{2}}{b_{0}^{2}}\left\| u\left( .,T\right) \right\|
_{L^{2}(0,l)}^{2}\mathbf{\leq }\frac{b_{1}^{2}b_{2}^{2}}{b_{0}^{2}}\left\|
u\right\| _{L^{2}(Q_{s})}^{2}\mathbf{+}\frac{b_{1}^{2}b_{2}^{2}}{b_{0}^{2}}%
\left\| u_{t}\right\| _{L^{2}(Q_{s})}^{2}\text{,}  \tag{28}
\end{equation}

we estimate the last term of the right-hand side of the inequality(27),we
obtain 
\begin{gather}
\left\| \Im _{x}^{*}u_{tt}(.,s)\right\| _{L^{2}(0,l)}^{2}+  \tag{29} \\
+\frac{b_{0}}{2}\left\| u_{t}(.,T)\right\| _{L^{2}(0,l)}^{2}+\beta \left\|
u_{tt}\left( .,s\right) \right\| _{L^{2}(0,l)}^{2}  \notag \\
\leq \left( \frac{2b_{1}^{2}+b_{2}}{b_{0}}\right) \left\| \Im
_{x}^{*}u_{tt}\right\| _{L^{2}(Q_{s})}^{2}+\frac{\beta b_{2}}{b_{0}}\left\|
u_{tt}\right\| _{L^{2}(Q_{s})}^{2}  \notag \\
+\frac{\frac{2b_{1}^{2}b_{2}^{2}}{b_{0}^{2}}+2b_{1}b_{2}+\left(
b_{2}+b_{1}\right) ^{2}+b_{3}^{2}}{b_{0}}\left\| u_{t}\right\|
_{L^{2}(Q_{s})}^{2}  \notag \\
+\frac{b_{0}^{2}\left( b_{2}^{2}+b_{4}^{2}+b_{5}^{2}\right)
+2b_{1}^{2}b_{2}^{2}}{b_{0}^{3}}\left\| u\right\| _{L^{2}(Q_{s})}^{2}  \notag
\end{gather}

For estimate the last term of the right-hand side of the inequality(29),we
will be proove the inequality $\left\| u\right\| _{L^{2}\left( Q_{s}\right)
}^{2}\mathbf{\leq }24T^{2}\left\| u_{t}\right\| _{L^{2}\left( Q_{s}\right)
}^{2},$ combining the last inequality and (29) we get

\begin{gather}
\left\| \Im _{x}^{*}u_{tt}(.,s)\right\| _{L^{2}(0,l)}^{2}+\left\|
u_{tt}\left( .,s\right) \right\| _{L^{2}(0,l)}^{2}+  \tag{30} \\
\;\;\;\;\;\;\;\;\;+\left\| u_{t}(.,T)\right\| _{L^{2}(0,l)}^{2}  \notag \\
\leq k\left[ \left\| \Im _{x}^{*}u_{tt}\right\| _{L^{2}(Q_{s})}^{2}+\left\|
u_{tt}\right\| _{L^{2}(Q_{s})}^{2}+\left\| u_{t}\right\| _{L^{2}(Q_{s})}^{2}%
\right] ,  \notag
\end{gather}

where$\;k\mathbf{=}\frac{\max \left( \beta b_{2},\left(
2b_{1}^{2}+b_{2}\right) ,b_{0}k\left( b_{i},T\right) \right) }{b_{0}\min
\left( 1,\beta ,\frac{b_{0}}{2}\right) }$.\newline
$k\left( b_{i},T\right) =$\ $\frac{2b_{1}^{2}b_{2}^{2}+b_{0}^{2}\left[
2b_{1}b_{2}+\left( b_{2}+b_{1}\right) ^{2}+b_{3}^{2}\right] +24T^{2}\left[
b_{0}^{2}\left( b_{2}^{2}+b_{4}^{2}+b_{5}^{2}\right) +2b_{1}^{2}b_{2}^{2}%
\right] }{b_{0}^{3}}$.\\
To continue,we introduce the new function $v(x,t)$\ $%
=\int_{t}^{T}u_{\tau \tau }d\tau $,\ then $u_{t}(x,t)=v(x,s)-v(x,t)$ , and $%
u_{t}(x,T)=v(x,s)$.\\
The inequality (30) it be 
\begin{gather}
\left\| \Im _{x}^{*}u_{tt}(.,s)\right\| _{L^{2}(0,l)}^{2}+\left\|
u_{tt}\left( .,s\right) \right\| _{L^{2}(0,l)}^{2}+  \tag{31} \\
+\left( 1-2k(T-s)\right) \left\| v(.,s)\right\| _{L^{2}(0,l)}^{2}  \notag \\
\leq 2k\left( \left\| \Im _{x}^{*}u_{tt}\right\| _{L^{2}(Q_{s})}^{2}+\left\|
u_{tt}\right\| _{L^{2}(Q_{s})}^{2}+\left\| v\right\|
_{L^{2}(Q_{s})}^{2}\right) \text{.}  \notag
\end{gather}
If $s_{0}>0$\ satisfies $\left( 1-2k(T-s_{0})\right) =\frac{1}{2}$,then the
inequality (31) implies 
\begin{gather}
\left\| \Im _{x}^{*}u_{tt}(.,s)\right\| _{L^{2}(0,l)}^{2}+\left\|
u_{tt}\left( .,s\right) \right\| _{L^{2}(0,l)}^{2}  \tag{32} \\
\;\;\;\;\;\;\;\;+\left\| v(.,s)\right\| _{L^{2}(0,l)}^{2}  \notag \\
\leq 4k\left( \left\| \Im _{x}^{*}u_{tt}\right\| _{L^{2}(Q_{s})}^{2}+\left\|
u_{tt}\right\| _{L^{2}(Q_{s})}^{2}+\left\| v\right\|
_{L^{2}(Q_{s})}^{2}\right) \text{,}  \notag
\end{gather}
for all $s\in \left[ T-s_{0},T\right]$.We denote

\vspace*{0.3cm}
\begin{center}
$Y(s)=\left\| \Im _{x}^{*}u_{tt}\right\| _{L^{2}(Q_{s})}^{2}+\left\| u_{tt}\right\|
_{L^{2}(Q_{s})}^{2}+\left\| v\right\| _{L^{2}(Q_{s})}^{2}$.
\end{center}
We get:\newline

\vspace*{0.2cm}
$Y^{\prime }(s)=-\left\| \Im _{x}^{*}u_{tt}(.,s)\right\|
_{L^{2}(0,l)}^{2}-\left\| u_{tt}\left( .,s\right) \right\|
_{L^{2}(0,l)}^{2}-\left\| v(.,s)\right\| _{L^{2}(0,l)}^{2}.$\newline

\newpage
Then and from (32) we obtain $-Y^{\prime }(s)\leq 4kY(s)\mathbf{.}$\\ 
Then $-\frac{\partial }{\partial s}\left( Y(s)\exp (4ks)\right) \leq 0$.\newline

\vspace*{0.3cm}
Integrating this inequality on $\left( s,T\right) $\ and taking into account
that $Y(T)=0$, we obtain $Y(s)\exp (4ks)\leq 0$.\\
Then $Y(s)=0$\ for all $s\in \left[ T-s_{0},T\right] .$Then $\omega =0$\ almost everywhere in 
$ Q_{T-s_{0}} $, proceding in this way step by step,we proove that $\omega =0\;
$ almost\ everywhere in $Q$.\newline

\vspace*{0.3cm}
This achieves the proof of proposition.Now we return to proove the th\'{e}or%
\`{e}me.We will be proove that $\overline{R(L)}=H$.\\
Since $H$\ is a Hilbert space ,the equality $\overline{R(L)}=H$\ is true,if from 
\begin{equation}
\left( Lu,W\right) _{H}=\left( \mathcal{L}u,\omega \right)
_{L^{2}(Q)}+\left( l_{1}u,\omega _{1}\right) _{L^{2}(0,l)}+\left(
l_{2}u,\omega _{2}\right) _{L^{2}(0,l)}=0,  \tag{33}
\end{equation}

where  $ W=\left( \omega ,\omega _{1},\omega _{2}\right) \in R(L)^{\perp }$,we get $\omega \equiv 0$,$\omega _{1}\equiv 0$\ and $\omega _{2}\equiv 0$\ in $Q$,for any element of $D_{0}(L) $.\newline
From (33) we obtain $\forall u\in D_{0}(L),\left( \mathcal{L}u,\omega
\right) _{L^{2}(Q)}=0$.Then by virtue of the 2$^{nd}$ proposition , we
conclude that $\omega \equiv 0$.\\
Then for (33),we obtain$\;\left(l_{1}u,\omega _{1}\right) _{L^{2}(0,l)}\mathbf{+}\left( l_{2}u,\omega
_{2}\right) _{L^{2}(0,l)}\mathbf{=}0$.\newline

\vspace*{0.3cm}
Since the quantities $l_{1}u$\ and $l_{2}u$\ can vanish independently and
the ranges of the trace operators $l_{1}$\ and $l_{2}$\ are dense in the
Hilbert space $L^{2}(0,l)$,then $\omega _{1}=\omega _{2}=0$\ .Thus to
conclude that $W=0$.

\end{document}